\documentclass[seceqn,secthm,draft]{elsart}
\journal{Insurance: Mathematics and Economics}
\usepackage{graphicx,amssymb,amsmath}
\usepackage{natbib}
\newcommand{\eps}{\varepsilon}
\newcommand{\floor}[1]{\lfloor{#1}\rfloor}
\newcommand{\ceil}[1]{\lceil{#1}\rceil}

\newcommand{\halmos}{\mbox{} \hfill $\Box$}
\newcommand{\1}{\mathbf{1}}

\newcommand{\E}{\mathrm{E}}

\newcommand{\be}{\begin{equation}}
\newcommand{\ee}{\end{equation}}

\newcommand{\Fbar}{\overline{F}}

\newenvironment{proof}{\paragraph*{\it Proof.}}{\halmos}

\begin{document}

\begin{frontmatter}

\title{Tails of random sums of a \\ heavy-tailed number of light-tailed terms\thanksref{ack}}

\thanks[ack]{Research supported financially by the AXA chair ``Large Risks in Insurance'' (C. Y. Robert) and by the IAP research network grant nr.\ P6/03 of the Belgian government (J. Segers). The authors gratefully acknowledge an anonymous referee for various suggestions leading to improved results throughout the paper.}

\author[rob]{Christian Y. Robert\corauthref{cor}}
\ead{chrobert@ensae.fr}
\corauth[cor]{Corresponding author.}
\address[rob]{ENSAE,  Timbre J120, 3 Avenue Pierre Larousse, 92245 MALAKOFF Cedex, France}

\author[seg]{Johan Segers}
\ead{Johan.Segers@uclouvain.be}
\address[seg]{Universit\'e catholique de Louvain, Institut de statistique, Voie du Roman Pays 20, B-1348 Louvain-la-Neuve, Belgium}

\begin{abstract}
The tail of the distribution of a sum of a random number of independent and identically distributed nonnegative random variables depends on the tails of the number of terms and of the terms themselves. This situation is of interest in the collective risk model, where the total claim size in a portfolio is the sum of a random number of claims. If the tail of the claim number is heavier than the tail of the claim sizes, then under certain conditions the tail of the total claim size does not change asymptotically if the individual claim sizes are replaced by their expectations. The conditions allow the claim number distribution to be of consistent variation or to be in the domain of attraction of a Gumbel distribution with a mean excess function that grows to infinity sufficiently fast. Moreover, the claim number is not necessarily required to be independent of the claim sizes.
\end{abstract}

\begin{keyword}
Compound distribution; Collective risk model; Consistent variation; Heavy-tailed distributions; Light-tailed distributions; Large deviations; Random sum.
\end{keyword}

\end{frontmatter}


\section{Introduction}

\label{S:intro}

In the collective risk model in actuarial risk theory, the total claim amount in a portfolio is assumed to be a random sum,
\[
    S_{N}=\sum_{i=1}^{N}X_{i}
\]
where $N$, the claim number, is a nonnegative integer-valued random
variable and $X_{1},X_{2},\ldots $, the claim sizes, form a sequence
of independent and identically distributed (iid) nonnegative random
variables, independent of $N$; see e.g. \citet[chapter~3]{KGDD}. The
tail of the compound distribution depends on the tails of the claim
number and claim size  distributions. If both distributions are
light-tailed, that is, if $\E[e^{tN}]<\infty $ and
$\E[e^{tX_{1}}]<\infty $ for some $t>0$, saddlepoint approximation
techniques can be used to analyze the tail of $S_{N}$ \citep{Beard,
Jensen}. If  the individual claim distribution is subexponential and
if the claim number distribution is light-tailed, then \be
    \Pr [S_{N}>x]\sim \E[N]\Pr [X_{1}>x], \qquad x\to \infty \label{E:heavyX}
\ee%
\citep[Theorem~A3.20]{EKM}; see \citet{Denisov} for one-sided versions of eq.~\eqref{E:heavyX} under milder conditions.

In this paper we concentrate on the converse case when the tail of $S_{N}$ is dominated by the tail of $N$. This setting is relevant, for instance, for earthquake insurance, featuring a potentially large number of bounded claims, or in the case of individual unobserved heterogeneity, as well as in queueing theory, see Section~\ref{S:Applications}. We try to answer to the following question raised in \citet{Schmidli}: When does
\be
\Pr [S_{N}>x]\sim \Pr [\E[X_{1}]N>x],\qquad x\to
\infty \label{E:heavyN}
\ee%
hold? We consider a slightly more general framework because we do not necessarily assume the claim number to be independent of the claim sizes.

Theorems~1.3 and 1.4 in \citet{Stam} shed light on the difference between the two approximations in eqs.~\eqref{E:heavyX} and \eqref{E:heavyN}: if the claim number and claim size distributions have finite expectations and regularly varying tails, that is, if there exist $\alpha >1$ and $\beta >1$ such that
\[
\lim_{x\to \infty }\frac{\Pr [N>xy]}{\Pr [N>x]}=y^{-\alpha
},\qquad \lim_{x\to \infty }\frac{\Pr [X_{1}>xy]}{\Pr
[X_{1}>x]}=y^{-\beta }
\]%
for $y>0$, then, provided $N$ and $\{X_{i}\}_{i\geq 1}$ are independent, $\alpha >\beta $ implies \eqref{E:heavyX} while $\beta >\alpha $ implies \eqref{E:heavyN}. In other words, which of the two approximations \eqref{E:heavyX} or \eqref{E:heavyN} is valid depends on which of the two distributions has the heavier tail. Approximation \eqref{E:heavyN} will therefore arise when the tail of $N$ is heavier than the tail of $X$.

The main contribution of this paper consists of the four theorems in Sections~\ref{S:results:CV} and \ref{S:results:Gumbel} providing different sets of sufficient conditions for eq.~\eqref{E:heavyN} to hold. The two theorems in Section~\ref{S:results:CV} concern the case where the claim number distribution is of consistent variation and the claim size distribution has a finite moment of order $r$ for some $r>1$. The two theorems in Section~\ref{S:results:Gumbel} treat the case where the claim number distribution is in the Gumbel domain of attraction and the moment generating function of the claim size distribution is finite in a neighbourhood of the origin. An important factor is the speed at which the mean excess function $\E[N-x\mid N>x]$ tends to infinity. In the special case where $N$ is a discretized Weibull random variable with tail function $\Pr [N>x]\sim \exp (-x^{\beta })$ as $x\to \infty $ and shape parameter $0< \beta <1$, the conditions for Theorems~\ref{T:Gumbel} and \ref{T:Gumbel2} translate into $\beta < 1/3$ and $\beta < 1/2$, respectively. The exponent $1/2$ marks the lower boundary on the speed of growth of the mean excess function for which eq.~\eqref{E:heavyN} can still be expected to hold; see Section~\ref{S:discussion}. In Theorems \ref{T:CV} and \ref{T:Gumbel}, the dependence between the claim number 
and the claim sizes 
can be arbitrary.

The outline of this paper is as follows: Section~\ref{S:prelimin} presents 
preliminaries on tail asymptotics and large deviations. Statements and proofs of our main results are to be found in Sections~\ref{S:results:CV} and \ref{S:results:Gumbel}. Section~\ref{S:discussion} gives a short discussion. Section~\ref{S:Applications} concludes the paper with applications to insurance and operations research.

The following notations and assumptions are in force throughout the
paper. Let $N$ be a nonnegative integer-valued random variable and
let $X_{1},X_{2},\ldots $ be a sequence of iid nonnegative random
variables with finite mean $\mu =\E[X_{1}]$ and variance $\sigma
^{2}=\mathrm{Var}[X_{1}]$. The claim number $N$ is not necessarily
assumed to be independent of $\{X_{i}\}_{i\geq 1}$. Denote $S_{0}=0$
and $S_{n}=X_{1}+\cdots +X_{n}$ for an integer $n\geq 1$. For real
$y$, we denote by $\lceil {y}\rceil $ the smallest integer greater
than or equal to $y$; similarly, $\lfloor {y}\rfloor $ is the
greatest integer smaller than or equal to $y$. For two positive
functions $f$ and $g$ defined in a neighbourhood of infinity, we
write $f(x)\sim g(x)$ as $x \to \infty $ provided $\lim_{x\to \infty
}f(x)/g(x)=1$ and we write $f(x)=o\left( g(x)\right) $ as $x\to
\infty $ provided $\lim_{x\to \infty }f(x)/g(x)=0$. The indicator
function of an event $A$ is denoted by $\1(A)$.


\section{Preliminaries}

\label{S:prelimin}


\subsection{Consistent variation}

\label{SS:CV}

A distribution function $F$ is said to have a \emph{consistently varying tail}, notation $F\in \mathcal{C}$, if
\be
\lim_{y\uparrow 1}\limsup_{x\to \infty }\frac{\Fbar(xy)}{%
\Fbar(x)}=1,  \label{E:CV}
\ee%
where $\Fbar=1-F$. This regularity property was first introduced in \citet{Cline} and
called ``intermediate regular variation''. The concept of consistent
variation has
been proven useful in various papers in queueing systems and ruin theory %
\citep{Cai, Kaas, Ng, Tang}.

A distribution function $F$ is said to have a \emph{regularly
varying tail}, notation $F\in \mathcal{R}$, if there exists $\alpha
>0$ such that
\[
\lim_{x\to \infty }\frac{\Fbar(xy)}{\Fbar(x)}%
=y^{-\alpha },\qquad y>0,
\]%
or equivalently such that
\[
\Fbar(x)=x^{-\alpha }L\left( x\right) ,\qquad x>0,
\]%
where $L$ is a slowly varying function, that is, a function satisfying%
\[
    \lim_{x\to \infty }\frac{L(xy)}{L(x)}=1,\qquad y>0.
\]
In \citet{ClineS}, it is shown that the class $%
\mathcal{C}$ strictly contains the class $\mathcal{R}$. In addition
the class $\mathcal{C}$ is itself strictly contained in the class of
distribution functions with a \emph{dominatedly varying tail}, that
is, the class $\mathcal{D}$ of distribution functions $F$ such that
\be \label{E:DV}
    \limsup_{x\to \infty }\frac{\Fbar(xy)}{\Fbar(x)} < \infty
\ee
for all (or, equivalently, for some) $y\in (0,1)$.

For a distribution function $F$, define
\[
    \Fbar_{\ast}(y) = \liminf_{x \to \infty} \frac{\Fbar(xy)}{\Fbar(x)}, \qquad y > 0,
\]
and
\[
    \alpha _{F} = \lim_{y \to \infty} \frac{- \log \{\Fbar_{\ast}(y)\}}{\log y} < \infty,
\]
with $- \log 0 = \infty$. In the terminology of \citet[Proposition~2.2.5]{BGT}, the quantity $\alpha_{F}$ is the \emph{upper Matuszewska index} of the function $1/\Fbar$. It is an upper bound for the \emph{upper order}, $\rho_{F}$, of $1/\Fbar$, defined by
\be
\label{E:rhoF}
    \rho_{F} = \limsup_{x \to \infty} \frac{-\log \{\Fbar(x)\}}{\log x}.
\ee
By \citet[Theorem~2.1.7]{BGT}, a distribution function $F$ belongs to $\mathcal{D}$ if and only if $\alpha_F < \infty$. But then $\rho_F \leq \alpha_F < \infty$. As $\Fbar(x) = x^{-\rho(x)}$ with $\limsup_{x \to \infty} \rho(x) = \rho_F$, we obtain that for every $\rho > \rho_F$,
\be
\label{E:rho}
    x^{-\rho } = o(\Fbar(x)), \qquad x \to \infty.
\ee
In this sense, distributions of dominated variation (and thus also those of consistent variation) are heavy-tailed.


\subsection{Domains of attraction}

\label{SS:DA}

Let $N$ be a nonnegative integer-valued random variable. Assume that
$N$ is in the maximal domain of attraction of an extreme value
distribution.
Because the distribution of $N$ is discrete, its attractor must be a Fr\'{e}%
chet or the Gumbel distribution, that is, its attractor cannot be a
Weibull extreme value distribution.

On the one hand, the distribution of $N$ is in the domain of
attraction of the Fr\'{e}chet distribution with shape parameter
$\alpha >0$ if and only if
\be
\Pr [N>x]=x^{-\alpha }L\left( x\right)  \label{E:RV}
\ee
where $L$ is a slowly varying function. This case was treated in
\citet{Stam} and in \citet{FGMS}; see Section~\ref{S:Applications}.
A sufficient
condition for \eqref{E:RV} is the von Mises type condition $%
\lim_{n\to \infty }n\Pr [N=n]/\Pr [N>n]=\alpha $
\citep{Anderson}.
In turn, regular variation implies consistent variation, that is, %
\eqref{E:RV} implies \eqref{E:CV} with $F(x)=\Pr [N\leq x]$, and it
is merely the latter concept which will be used later on.

On the other hand, the distribution of $N$ is in the domain of
attraction of the Gumbel distribution if and only if there exists a
positive function $a$ defined in a neighbourhood of infinity such
that
\be
\lim_{x\to \infty }\frac{\Pr [N>x+ya(x)]}{\Pr
[N>x]}=e^{-y},\qquad y\in \mathbb{R}.  \label{E:Gumbel}
\ee%
The function $a$ is necessarily asymptotically equivalent to the
mean excess function, $a(x)\sim \E[N-x\mid N>x]$ as
$x\to \infty $, and it is self-neglecting, that is,
\be
\lim_{x\to \infty }\frac{a(x)}{x}=0,\qquad
\lim_{x\to \infty }\frac{a(x+ya(x))}{a(x)}=1,\qquad y\in
\mathbb{R}.  \label{E:SelfNeglecting}
\ee%
The tail function of $N$ admits the following useful representation:
\be
\Pr [N>x]=c(x)\exp \biggl(-\int_{0}^{x}\frac{g(t)}{a(t)}dt\biggr)
\label{E:repr}
\ee%
where $\lim_{x\to \infty }c(x)=c\in (0,\infty )$ and $%
\lim_{t\to \infty }g(t)=1$. By \citet{Anderson}, a
sufficient
condition for $N$ to be in the Gumbel domain of attraction is that $%
q_{n}\to \infty $ and $q_{n+1}-q_{n}\to 0$ as
$n\to \infty $, where $q_{n}=\Pr [N>n]/\Pr [N=n]$.

\paragraph*{Example.}
Let $N$ be a discretized Weibull variable, that is, $N=\lfloor
{Y}\rfloor $
where $\Pr [Y>y]=\exp (-y^{\beta })$ for $y\geq 0$, with parameter $\beta >0$%
. If $\beta <1$, then $\Pr [N>y]\sim \Pr [Y>y]$ as $y\to
\infty $
and \eqref{E:Gumbel}--\eqref{E:repr} hold with $a(x)=\beta ^{-1}x^{1-\beta }$%
.


\subsection{Large deviations}

\label{SS:LD}

In order to analyse the tails of the compound sum $S_{N}$, we need
bounds on the probability of large deviations of $S_{n}$. The first
bound is a special case of Lemma~2.3 in \citet{Tang2}.

\begin{lem}
\label{L:Tang} If $\E[X_{1}^{r}]<\infty $ for some $r>1$,
then for each $\gamma >0$ and $q>0$, there exist positive numbers
$v$ and $C=C\left( v,\gamma \right) $ irrespective to $x$ and $n$
such that for all $x\geq \gamma n$ and $n=1,2,\ldots ,$
\[
\Pr [S_{n}-n\mu >x]\leq n\Pr \left[ X_{1}-\mu >vx\right] +Cx^{-q}.
\]
\end{lem}

The second result is a simple consequence of Cram\'{e}r's theorem on
large deviations; see \citet[equation~(7.30), p.~553]{Feller} or \citet[Lemma~2.2]{BDK}.

\begin{lem}
\label{T:LD} If $\E[e^{tX_{1}}]<\infty $ for some $t>0$,
then for
any sequence $a_{n}$ satisfying $a_{n}/n^{1/2}\to \infty $ and $%
a_{n}/n\to 0$ as $n\to \infty $,
\[
\begin{array}{rcl}
\Pr [S_{n}>n\mu +a_{n}] & = & \displaystyle\exp \biggl(-\frac{1}{2\sigma ^{2}%
}\frac{a_{n}^{2}}{n}\{1+o(1)\}\biggr), \\[1em]
\Pr [S_{n}<n\mu -a_{n}] & = & \displaystyle\exp \biggl(-\frac{1}{2\sigma ^{2}%
}\frac{a_{n}^{2}}{n}\{1+o(1)\}\biggr).%
\end{array}%
\]
\end{lem}

\section{Main results: Consistent variation}
\label{S:results:CV}

In this section, we treat the case where the tail of the distribution function of $N$ is of consistent variation, see eq.~\eqref{E:CV}. Theorem~\ref{T:CV} states a simple sufficient condition on the common distribution of the $X_i$ for eq.~\eqref{E:heavyN} to be valid for arbitrary dependence structure between $N$ and $\{X_i\}$. In case of independence, a weaker condition suffices, see Theorem~\ref{T:CV2}. The results crucially rest on a large-deviations result by \citet[Lemma~2.3]{Tang2}, reproduced as Lemma~\ref{L:Tang} above, and which was kindly pointed out to us by an anonymous Referee.

\begin{thm}
\label{T:CV}
If the tail of the distribution of $N$ is of consistent variation, if $\E[X_1^r] < \infty$ for some $r > 1$ and if
\be
\label{E:p}
    x \Pr [X_{1} > x] = o( \Pr [N > x] ), \qquad x \to \infty,
\ee
then $\Pr [S_{N}>x] \sim \Pr [N>x/\mu ]$ as $x \to \infty$.
\end{thm}

\begin{proof}
Since $\Pr [S_{N}>x]=\Pr [S_{N}/\mu >x/\mu ]$, we can without loss of generality assume that $\mu =1$. Fix $0<\eps <1$. For $x>0$,
\be
\label{E:CV:5}
    \Pr [S_{N}>x] = \Pr [S_{N}>x, \, N\leq (1-\eps )x] + \Pr[S_{N}>x, \, N>(1-\eps)x],
\ee so that, on the one hand,

\be \label{E:CV:10}
    \Pr[S_N > x] \leq \Pr [S_{\lfloor {(1-\eps )x}\rfloor }>x] + \Pr[N > (1-\eps) x],
\ee
and, on the other hand,
\begin{eqnarray}
\label{E:CV:20}
    \Pr[S_N > x]
    &\geq& 0 + \Pr [S_{\lceil {(1+\eps )x} \rceil } > x, \, N > (1+\eps )x] \nonumber \\
    &\geq& \Pr [N > (1+\eps )x] - \Pr [S_{\lceil {(1+\eps)x}\rceil} \leq x].
\end{eqnarray}
Since the tail of $N$ is of consistent variation, by
eq.~\eqref{E:CV}, \be \label{E:CV:30}
    \lim_{\eps \downarrow 0} \liminf_{x \to \infty} \frac{\Pr [N>(1+\eps )x]} {\Pr [N>(1-\eps )x]} = 1.
\ee
In view of eqs.~\eqref{E:CV:10}, \eqref{E:CV:20} and \eqref{E:CV:30}, $\Pr[S_N > x] \sim \Pr[N > x]$ as $x \to \infty$ will follow if we can show that
\begin{eqnarray}
\label{E:CV:40a}
    \Pr [S_{\lfloor {(1-\eps )x}\rfloor} > x  ] &=& o(\Pr[N > x]), \qquad x \to \infty; \\
\label{E:CV:40b}
    \Pr [S_{\lceil {(1+\eps)x}  \rceil} \leq x] &=& o(\Pr[N > x]), \qquad x \to \infty.
\end{eqnarray}

We first show eq.~\eqref{E:CV:40a}. Let $\alpha_{N}$ be the upper Matuszewska index of the function $x\mapsto 1/\Pr [N>x]$. Since the tail of $N$ is of consistent variation, $\alpha_N < \infty$ (see subsection~\ref{SS:CV}). Pick $\rho > \alpha_N$. By eq.~\eqref{E:rhoF},
\be
\label{E:CV:50}
    x^{-\rho} = o(\Pr[N > x]), \qquad x \to \infty.
\ee
By Lemma~\ref{L:Tang}, there exist positive numbers $v$ and $C$ such that
\begin{eqnarray*}
    \Pr [ S_{\lfloor {(1-\eps )x}\rfloor }>x]
    &\leq& \Pr [ S_{\lfloor {(1-\eps )x}\rfloor }-\lfloor {(1-\eps )x} \rfloor > \eps x]  \\
    &\leq& \lfloor {(1-\eps )x}\rfloor \Pr [X_{1} > v \eps x] + C(\eps x)^{-\rho}
\end{eqnarray*}
and thus
\[
    \frac{\Pr [ S_{\lfloor {(1-\eps )x}\rfloor }>x]}{\Pr[N > x]}
    \leq \frac{1}{v{\eps }}
        \frac{v{\eps }x\Pr [X_{1}>v{\eps }x]}{\Pr [N>v{\eps }x]}
        \frac{\Pr [N>v{\eps }x]}{\Pr [N>x]}
        + \frac{C(\eps x)^{-\rho}}{\Pr[N > x]}.
\]
Eq.~\eqref{E:CV:40a} now follows from the above inequality combined with eqs.~\eqref{E:DV}, \eqref{E:p}, and \eqref{E:CV:50}.

Next we show eq.~\eqref{E:CV:40b}. Since the $X_i$ are nonnegative, by Chernoff's bound, there exists $0 < a < 1$, depending on $\eps > 0$, such that for all sufficiently large $x$,
\be
\label{E:CV:60}
    \Pr[ S_{\ceil{(1+\eps)x}} \leq x ] \leq a^x.
\ee
This inequality in combination with \eqref{E:CV:50} yields eq.~\eqref{E:CV:40b}.
\end{proof}

\begin{thm}
\label{T:CV2}
Assume that $N$ and $\{X_i\}_{i \geq 1}$ are independent. If the tail of the distribution of $N$ is of consistent variation, if $\E[X_1^r] < \infty$ for some $r > 1$ and if one of the following two conditions holds:
\begin{description}
\item[{\it case $\E[N] < \infty$:}]
\be
\label{E:CV2:finite}
    \Pr [ X_1 > x ] = o(\Pr [N > x]), x \to \infty;
\ee
\item[{\it case $\E[N] = \infty$:}] there exists $q$ with $1 \leq q < r$ such that
\be
\label{E:CV2:infinite}
    \limsup_{x\to \infty } \frac{\E[ N \1(N \leq x)]}{x^q \Pr[N > x]} < \infty;
\ee
\end{description}
then $\Pr [S_{N} > x] \sim \Pr [N>x/\mu ]$ as $x\to \infty$.
\end{thm}

\begin{proof}
The proof is similar to that of Theorem~\ref{T:CV}. We just indicate
the modifications. In view of eq.~\eqref{E:CV:5}, it suffices to
show that \be \label{E:CV2:10}
    \Pr [ S_{N}>x, \, N \leq (1-\eps )x] = o(\Pr[N > x]), \qquad x \to \infty.
\ee
Since $N$ and $\{X_{i}\}_{i\geq 1}$ are independent,
\[
    \Pr [ S_{N}>x, \, N\leq (1-\eps )x]
    = \sum_{k=0}^{\lfloor (1-\eps )x \rfloor}\Pr [N = k] \Pr [S_{k} - k > x - k].
\]
Pick $\rho > \alpha_N$, the upper Matuszewska index of the function $x\mapsto 1/\Pr [N>x]$. By Lemma~\ref{L:Tang}, there exist positive numbers $v$ and $C$ such that
\begin{eqnarray*}
    \lefteqn{
    \Pr [ S_{N}>x, \, N \leq (1-\eps )x ]
    } \nonumber \\
    &\leq &\sum_{k=0}^{\lfloor {(1-\eps )x}\rfloor }\Pr [ N=k]
    \biggl( k \Pr [ X_{1}>v( x-k) ]+  \frac{C}{(x - k)^\rho} \biggr) \nonumber \\
    &\leq & \Pr [ X_{1}>\eps vx] \sum_{k=0}^{\lfloor (1-\eps )x\rfloor }\Pr [ N = k ] k
    + \frac{C}{(\eps x)^\rho} \nonumber \\
    &=& \Pr [ X_{1} > \eps vx] \E [ N \1\{N \leq (1-\eps)x\} ] + \frac{C}{(\eps x)^\rho}.
\end{eqnarray*}
In view of eq.~\eqref{E:CV:50}, we only need to deal with the first term on the right-hand side of the previous display. We now invoke the additional condition. On the one hand, if $\E[N] < \infty$, then also $\E [ N \1\{N \leq (1-\eps)x\} ] \leq \E[N] < \infty$, and
\[
    \frac{\Pr [ X_{1}>\eps vx]}{\Pr[N > x]}
    = \frac{\Pr[N > \eps v x]}{\Pr[N > x]} \frac{\Pr[X_1 > \eps v x]}{\Pr[N > \eps v x]}
    \to 0, \qquad x \to \infty,
\]
where we used eqs.~\eqref{E:DV} and \eqref{E:CV2:finite}. On the other hand, if $\E[N] = \infty$, then, as $\Pr [X_1 > x] \leq \E[X_1^r] / x^r$,
\begin{eqnarray*}
    \Pr[X_1 > \eps vx] \E [ N \1\{N \leq (1-\eps)x\} ]
    &=& O(x^{-r}) O(x^q \Pr[N > (1-\eps)x]) \\
    &=& o(\Pr[N > x]), \qquad x \to \infty,
\end{eqnarray*}
as required.
\end{proof}

\paragraph*{Remark.}
Observe that if $\E[N] = \infty$ and the function $x \mapsto \Pr[N > x]$ is regularly varying of index $-\alpha$ for some $\alpha \in (0, 1]$, then eq.~\eqref{E:CV2:infinite} holds, for, by Karamata's theorem \citep[Propositions~1.5.8 and 1.5.9a]{BGT},
\[
    \E[N \1(N \leq x)] \leq \int_0^x \Pr[N > y] dy
    \left\lbrace
    \begin{array}{l@{\qquad}l}
    \sim x \Pr[N > x] / (1-\alpha) & \mbox{if $0 < \alpha < 1$,} \\
    = o(x^q \Pr[N > x]) & \mbox{if $\alpha = 1 < q$.}
    \end{array}
    \right.
\]

\section{Main results: Gumbel domain of attraction}
\label{S:results:Gumbel}

Next, we treat the case where the claim number distribution is in the Gumbel domain of attraction (see subsection~\ref{SS:DA}) and the moment generating function of the claim size distribution is finite in a neighbourhood of the origin. An important factor is the speed at which the auxiliary function $a$ in \eqref{E:Gumbel} tends to infinity. If $a(x)/x^{2/3}\to \infty $ as $x\to \infty $, then \eqref{E:heavyN} holds without further conditions (Theorem~\ref{T:Gumbel}). If the function $a$ is merely assumed to have a \emph{lower order} larger than $1/2$ in the terminology of \citet[Section~2.2.2, p.~73]{BGT}, then \eqref{E:heavyN} still holds provided the claim number is independent of the claim sizes (Theorem~\ref{T:Gumbel2}). The order $1/2$ marks the lower boundary on the speed of growth of $a$ for which \eqref{E:heavyN} can still be expected to hold true; see Section~\ref{S:discussion}.

\begin{thm}
\label{T:Gumbel} If $\E[e^{\gamma X_{1}}]<\infty $ for some
$\gamma >0$ and if \eqref{E:Gumbel} holds for a function $a$ such
that
\be
\frac{a(x)}{x^{2/3}}\to \infty ,\qquad x\to \infty ,
\label{hypa1}
\ee%
then $\Pr [S_{N}>x]\sim \Pr [N>x/\mu ]$ as $x\to \infty $.
\end{thm}

\begin{proof}
Without loss of generality, assume $\mu =1$. Fix $\eps >0$.
Write
\be
\Pr [S_{N}>x]=\Pr [S_{N}>x,N\leq x-\eps a(x)]+\Pr
[S_{N}>x,N>x-\eps a(x)].  \label{dec2}
\ee%
By Lemma~\ref{T:LD}, the first term on the right-hand side of
\eqref{dec2} is bounded by
\[
\begin{array}{rcl}
\Pr [S_{N}>x,N\leq x-\eps a(x)] & \leq & \Pr [S_{\lfloor {%
x-\eps a(x)}\rfloor }>x] \\[1em]
& = & \displaystyle\exp \biggl(-\frac{\eps ^{2}}{2\sigma ^{2}}\frac{%
a^{2}(x)}{x}\{1+o(1)\}\biggr),\qquad x\to \infty .%
\end{array}%
\]%
By \eqref{hypa1}, $a^{2}(x)/x$ is of larger order than $%
\int_{0}^{x}a^{-1}(t)dt$, so that the right-hand side in the
previous display must be $o(\Pr [N>x])$ as $x\to \infty $.

The second term on the right-hand side of \eqref{dec2} is bounded
from above by
\[
\Pr [S_{N}>x,N>x-\eps a(x)]\leq \Pr [N>x-\eps a(x)]
\]%
and from below by
\[
\begin{array}{rcl}
\Pr [S_{N}>x,N>x-\eps a(x)] & \geq & \Pr [S_{\lceil
{x+\eps
a(x)}\rceil }>x,N>x+\eps a(x)] \\[1ex]
& \geq & \Pr [N>x+\eps a(x)]-\Pr [S_{\lceil {x+\eps a(x)}%
\rceil }\leq x].%
\end{array}%
\]%
The second term on the right-hand side is $o(\Pr [N>x])$ as
$x\to
\infty $ by the same argument as in the previous paragraph. Moreover, by %
\eqref{E:Gumbel},
\[
\lim_{x\to \infty }\frac{\Pr [N>x+\eps a(x)]}{\Pr
[N>x-\eps a(x)]}=e^{-2\eps }.
\]%
As $\eps >0$ was arbitrary, indeed $\Pr [S_{N}>x]\sim \Pr [N>x]$ as $%
x\to \infty $.
\end{proof}

\begin{thm}
\label{T:Gumbel2} If $N$ and $\{X_{i}\}_{i\geq 1}$ are independent, if $\E%
[e^{\gamma X_{1}}]<\infty $ for some $\gamma >0$ and if
\eqref{E:Gumbel} holds for a function $a$ such that
\be
\liminf_{x\to \infty }\frac{\log a(x)}{\log x}>\frac{1}{2},
\label{E:a}
\ee%
then $\Pr [S_{N}>x]\sim \Pr [N>x/\mu ]$ as $x\to \infty $.
\end{thm}

\begin{proof}
Without loss of generality, assume $\mu =1$. Let $\lambda (x):=\log
a(x)/\log x$, so $a(x)=x^{\lambda (x)}$. By \eqref{E:SelfNeglecting} and %
\eqref{E:a},
\[
1/2<\lambda _{0}:=\liminf_{x\to \infty }\lambda (x)\leq
\limsup_{x\to \infty }\lambda (x)\leq 1.
\]%
Let $k$ be an integer larger than $2$ such that $\lambda
_{0}>1/2+1/2^{k} $ and decompose \be \Pr
[S_{N}>x]=T_{1}(x)+T_{2}(x)+T_{3}(x),  \label{E:ST}
\ee%
where
\[
\begin{array}{rcl}
T_{1}(x) & = & \Pr [S_{N}>x,N\leq \lfloor {x-x^{3/4}}\rfloor ], \\[1ex]
T_{2}(x) & = & \Pr [S_{N}>x,\lfloor {x-x^{3/4}}\rfloor <N\leq \lfloor {%
x-x^{1/2+1/2^{k}}}\rfloor ], \\[1ex]
T_{3}(x) & = & \Pr [S_{N}>x,\lfloor {x-x^{1/2+1/2^{k}}}\rfloor <N].%
\end{array}%
\]%
We will show that for $i=1,2$,
\be
T_{i}(x)=o(\Pr [N>x]),\qquad x\to \infty ,  \label{E:T12}
\ee%
and, for arbitrary $\eps >0$,
\be
e^{-\eps }\leq \liminf_{x\to \infty }\frac{T_{3}(x)}{\Pr [N>x]%
}\leq \limsup_{x\to \infty }\frac{T_{3}(x)}{\Pr [N>x]}\leq
e^{\eps }.  \label{E:T3}
\ee%
Since $\eps $ is arbitrary, the combination of equations \eqref{E:ST}%
, \eqref{E:T12} and \eqref{E:T3} implies $\lim_{x\to \infty
}\Pr [S_{N}>x]/\Pr [N>x]=1$, as required.

\textit{The term $T_1(x)$.}

Denoting $n(x)=\lfloor {x-x^{3/4}}\rfloor $, we have $T_{1}(x)\leq
\Pr [S_{n(x)}>x]=\Pr [S_{n(x)}>n(x)+\{x-n(x)\}]$. Since $x-n(x)\sim
x^{3/4}$ is
of smaller order than $x$ but of larger order than $x^{1/2}$ as $%
x\to \infty $, an application of Lemma~\ref{T:LD} yields
\[
T_{1}(x)\leq \exp \biggl(-\frac{1}{2\sigma ^{2}}x^{1/2}\{1+o(1)\}\biggr)%
,\qquad x\to \infty .
\]%
On the other hand, since $x^{1/2}=o(a(x))$ as $x\to \infty
$, the representation in \eqref{E:repr} implies $\exp (-\delta
x^{1/2})=o(\Pr
[N>x]) $ as $x\to \infty $ for all $\delta >0$. Equation~%
\eqref{E:T12} for $i=1$ follows.

\textit{The term $T_2(x)$.}

For integer $j\geq 2$, write $\tau _{j}=1/2+1/2^{j}$. Clearly
$3/4=\tau
_{2}>\tau _{3}>\cdots >\tau _{k}>1/2$. Further, denote $n_{j}(x)=\lfloor {%
x-x^{\tau _{j}}}\rfloor $. We have $n_{2}(x)\leq n_{3}(x)\leq \cdots
$ and
\be
\begin{array}[b]{rcl}
T_{2}(x) & = & \Pr [S_{N}>x,n_{2}(x)<N\leq n_{k}(x)] \\[1ex]
& = & \displaystyle\sum_{j=2}^{k-1}\Pr [S_{N}>x,n_{j}(x)<N\leq n_{j+1}(x)] \\%
[1em]
& \leq & \displaystyle\sum_{j=2}^{k-1}\Pr [S_{n_{j+1}(x)}>x]\Pr [N>n_{j}(x)].%
\end{array}
\label{E:T2:bound}
\ee%
A similar argument as for the term $T_{1}(x)$ in the previous
paragraph yields
\[
\begin{array}{rcl}
\Pr [S_{n_{j+1}(x)}>x] & = & \displaystyle\exp \biggl(-\frac{1}{2\sigma ^{2}}%
x^{2\tau _{j+1}-1}\{1+o(1)\}\biggr) \\[1em]
& = & \displaystyle\exp \biggl(-\frac{1}{2\sigma ^{2}}x^{1/2^{j}}\{1+o(1)\}%
\biggr),\qquad x\to \infty .%
\end{array}%
\]%
On the other hand, since $x^{1/2}=o(a(x))$ as $x\to \infty
$, the representation in \eqref{E:repr} implies
\[
\log \frac{\Pr [N>n_{j}(x)]}{\Pr [N>x]}=o(x^{1/2^{j}}),\qquad
x\to \infty .
\]%
Combine the final two displays to derive that every term on the
right-hand side of \eqref{E:T2:bound} is $o(\Pr [N>x])$ as
$x\to \infty $.

\textit{The term $T_3(x)$.}

Fix $\eps >0$. On the one hand, since $a(x)=x^{\lambda (x)}$,
the choice of $k$ entails that $x^{1/2+1/2^{k}}=o(a(x))$ as
$x\to \infty $. Hence
\[
T_{3}(x)=\Pr [S_{N}>x,N>\lfloor {x-x^{1/2+1/2^{k}}}\rfloor ]\leq \Pr
[N>x-\eps a(x)],
\]%
whence, by \eqref{E:Gumbel}, $T_{3}(x)\leq \{1+o(1)\}e^{\eps
}\Pr [N>x]$ as $x\to \infty $.

On the other hand, denoting $m(x)=\lfloor {x+\eps
a(x)}\rfloor $,
\be
\begin{array}[b]{rcl}
T_{3}(x) & \geq & \Pr [S_{N}>x,N>m(x)] \\[1ex]
& \geq & \Pr [S_{m(x)}>x,N>m(x)]=\{1-\Pr [S_{m(x)}\leq x]\}\Pr [N>m(x)].%
\end{array}
\label{E:T3:10}
\ee%
Since $m(x)=x+\{\eps +o(1)\}a(x)$ as $x\to \infty $, by Lemma~%
\ref{T:LD},
\be
\begin{array}[b]{rcl}
\Pr [S_{m(x)}\leq x] & = & \Pr [S_{m(x)}\leq m(x)-\{m(x)-x\}] \\[1em]
& = & \displaystyle\exp \biggl(-\frac{\eps ^{2}}{2\sigma ^{2}}\frac{%
a^{2}(x)}{x}\{1+o(1)\}\biggr)\to 0,\qquad x\to \infty .%
\end{array}
\label{E:T3:20}
\ee%
Moroever, by \eqref{E:Gumbel}, $\Pr [N>m(x)]\sim e^{-\eps
}\Pr [N>x]$ as $x\to \infty $. This relation in combination
with \eqref{E:T3:10} and \eqref{E:T3:20} yields $T_{3}(x)\geq
\{1+o(1)\}e^{-\eps }\Pr [N>x] $ as $x\to \infty $.
This finishes the proof of \eqref{E:T3} and hence of the theorem.
\end{proof}

\paragraph*{Example.}

Let $N$ be a discretized Weibull variable with shape parameter
$0<\beta <1$ as in the example in subsection~\ref{SS:DA}. Then
Theorem~\ref{T:Gumbel} applies for $\beta <1/3$, while Theorem~\ref{T:Gumbel2}
applies as long as $\beta <1/2$.


\section{Discussion}

\label{S:discussion}

At first sight, Theorem~\ref{T:Gumbel2} does not seem to entail much of an
extension compared to the following corollary to Theorem~3.6 in %
\citet{Asmussen}.

\begin{cor}
\label{C:Asmussen}
If $N$ and $\{X_{i}\}_{i\geq 1}$ are independent and if eq.~\eqref{E:Gumbel} holds for a function $a$ such that
\begin{description}
\item[{\it (i)}] $a(n) / n^{1/2} \to \infty$ as $n \to \infty$,
\item[{\it (ii)}] $\Pr [S_{n}>n\mu +ca(n)]=o(\Pr [N\geq n])$ as $%
n\to \infty $ for all $c>0$,
\end{description}
then $\Pr [S_{N}>x]\sim \Pr [N>x/\mu ]$ as $x\to \infty $.
\end{cor}

Indeed, the only real difference of Corollary~\ref{C:Asmussen} with respect to Theorem~\ref{T:Gumbel2} seems to be the extra condition (ii). However, this condition turns out to be not so harmless: Although condition~(i) only requires $a(n)$ to grow to infinity at a faster rate than $n^{1/2}$, the following lemma shows that condition~(ii) effectively forces a much faster rate on $a$, comparable to the one imposed in Theorem~\ref{T:Gumbel}.

\begin{lem}
\label{L:Asmussen}
Under the conditions of Corollary~\ref{C:Asmussen}, if $\E[e^{\gamma X_{1}}]<\infty $ for some $\gamma >0$ and if the function $a$ is regularly varying of index $\delta <1$, then necessarily $\delta \geq 2/3$.
\end{lem}

\begin{proof}
By \eqref{E:repr} and by Lemma~\ref{T:LD},
\be
\frac{\Pr [S_{n}>n\mu +a(n)]}{\Pr [N\geq n]}=\exp
\biggl(-\frac{1}{2\sigma
^{2}}\frac{a^{2}(n)}{n}\{1+o(1)\}+\int_{0}^{n}\frac{g(t)}{a(t)}dt-\log c(n)%
\biggr)  \label{E:Asmussen}
\ee%
as $n\to \infty $. By Karamata's theorem %
\citep[Proposition~1.5.8]{BGT},
\[
\int_{0}^{n}\frac{g(t)}{a(t)}dt\sim \frac{1}{1-\delta
}\frac{n}{a(n)},\qquad n\to \infty .
\]%
On the one hand, the function $t\mapsto a^{2}(t)/t$ is regularly
varying of index $2\delta -1$; on the other hand, the function
$t\mapsto t/a(t)$ is
regularly varying of index $1-\delta $. Hence, if the expression in %
\eqref{E:Asmussen} converges to zero as $n\to \infty $, then
necessarily $2\delta -1\geq 1-\delta $, whence $\delta \geq 2/3$.
\end{proof}

Finally, the lower bound $1/2$ for the order of the mean excess
function $a$ in Theorem~\ref{T:Gumbel2} seems to mark the minimal weight
that must be present in the tail of $N$ for the asymptotic
equivalence $\Pr [S_{N}>x]\sim \Pr [N>x/\mu ]$ to be true. Assume
for example that the distribution of $X_{i}$
is unit-mean exponential and that $N$ is independent of $\{X_{i}\}_{i\geq 1}$%
. For $t\geq 0$, the distribution of the random variable $Z_{t}=\max
\{n=0,1,\ldots :S_{n}\leq t\}$ is Poisson with mean $t$, and
\[
\Pr [S_{N}>t]=\Pr [N>Z_{t}]=\E[\exp \{-g(Z_{t})\}],
\]%
where $g(x):=-\log \Pr [N>x]$. The asymptotic behavior of the final
expression in the previous display as $t\to \infty $ has been
studied in \citet{Foss} in the general case where $Z_{t}$ is the sum
of iid nonnegative random variables; in our case, $Z_{t}$ is for
integer $t$ the sum of independent random variables with common
Poisson distribution and mean $1$. By \citet[Theorem 5.1]{Foss}, if
$N=\lfloor {Y}\rfloor $ and $\Pr [Y>y]=\exp (-y^{\beta })$ for
$y\geq 0$ where $\beta \in \lbrack 1/2,2/3)$, then
\[
\Pr [S_{N}>x]\sim \Pr [N>x]\exp (\beta ^{2}x^{2\beta -1}/2),\qquad
x\to \infty .
\]%
In particular, $\Pr [S_{N}>x]$ is not asymptotically equivalent to
$\Pr [N>x] $. In \citet{Asmussen}, the exponent $1/2$ was found to
be critical as well.


\section{Applications}
\label{S:Applications}

\paragraph*{Earthquake insurance.}
Earthquake insurance provides coverage to the policyholder in the
event of an earthquake that causes damage to the policyholder's
properties. Insurance companies must be careful when assigning this
type of insurance because, even if the individual claims are bounded
by the value of the properties insured, the number of claims can be
very large. An earthquake strong enough to destroy one house will
probably destroy hundreds of houses in the same area.

Let us assume that, given the energy of the earthquake
$\Lambda=\lambda $, the number of claims has a Poisson distribution
with parameter $\beta \lambda $ where $\beta $ is a positive
constant. Despite of the apparent complexity involved in the
dynamics of earthquakes, the probability distribution of the energy
of an earthquake follows a power law distribution known as the
Gutenberg-Richter law \citep{GuRi}: $\Pr [ \Lambda >\lambda ] =
\lambda ^{-\alpha } L(\lambda)$ where the exponent $\alpha $ is an
universal exponent close to $1$, universal in the sense that it does
not depend on a particular geographic area, and $L$ is a slowly
varying function.

We claim that as $x \to \infty$,
\[
    \Pr[N > x] \sim \Pr[\Lambda > x/\beta] = (x/\beta)^{-\alpha} L(x/\beta).
\]

The proof goes as follows. Without loss of generality, assume $\beta
= 1$. Let $F_\lambda$ denote the distribution function of the
Poisson distribution with mean $\lambda$. Note that if $0 < \lambda
< \mu < \infty$, then $1 - F_\lambda < 1 - F_\mu = 1 - F_\lambda
\ast F_{\mu-\lambda}$. Also let $Z_1, Z_2, \ldots$ denote iid random
variables with common Poisson distribution and mean $1$. Then
\begin{eqnarray*}
    \Pr[N > x]
    &=& \E[\Pr(N > x \mid \Lambda)] = \E[1 - F_\Lambda(x)] \\
    &\leq& \E[1 - F_{\ceil{\Lambda}}(x)] \\
    &=& \Pr \left[ \sum_{i=1}^{\ceil{\Lambda}} Z_i > x \right]
\end{eqnarray*}
and similarly
\[
    \Pr[N > x] \geq \Pr \left[ \sum_{i=1}^{\floor{\Lambda}} Z_i > x \right].
\]
Since the tail function of $\Lambda$ is regularly varying and since
$\Lambda - 1 < \floor{\Lambda} \leq \ceil{\Lambda} < \Lambda+1$, we
find as $\lambda \to \infty$
\[
    \Pr[\floor{\Lambda} > \lambda] \sim \Pr[\Lambda > \lambda] \sim \Pr[\ceil{\Lambda} >
    \lambda].
\]
The claim now follows from our Theorem~\ref{T:CV} and the above
upper and lower bounds on $\Pr[N > x]$.

By Theorem~\ref{T:CV}, the distribution of the total claim amount in
the portfolio has a regularly varying tail as well. Let us note that
it is not necessary to assume that the individual claim amounts are
independent of the claim number. In fact the energy of the
earthquake may also have an impact on the distribution of the
individual claim amounts.

\paragraph*{Hierarchical unobserved heterogeneity.}
Let us consider an hierarchical heterogeneity model for the number, $N$, of claims of a policyholder. First assume that, given $\Lambda = \lambda$, the claim number $N$ has a Poisson distribution with parameter $\lambda$. Secondly assume that,
given $V=v$, $\Lambda $ has an exponential distribution with mean $v>0$. It follows that
\[
    \Pr [N=n \mid V=v] = \frac{v}{1+v} \frac{1}{(1+v)^{n}},
\]
that is, given $V=v$, $N$ follows a Geometric distribution with succes probability parameter $v/(1+v)$. Thirdly assume that $\log(1+V)$ follows a Gamma distribution with shape parameter $\gamma > 0$ and scale parameter $c>0$. We deduce that
\[
    \Pr[ N \geq n] = \E[e^{-n\log(1+V)}] = \biggl( \frac{c}{c+n} \biggr) ^{\gamma },
\]%
that is, $N$ has a Pareto distribution with index $\gamma $.

It is well-known that the omission of an individual unobserved heterogeneity leads to overdispersion (in the sense that the variance is larger with heterogeneity). The proposed hierarchical model shows that unobserved heterogeneity can lead to grossly incorrect conclusions about the tail of the claim number distribution. Assume that their is no heterogeneity for the claim distributions and that the common distribution is light-tailed. If the heterogeneity for the number of claims is observed, then the distribution of the total claim amount is light-tailed, whereas if the heterogeneity is unobserved, then the distribution of the total claim amount is heavy-tailed.

\paragraph*{Stationary waiting time of customers.}
Let $\{T_{i}\}_{i\geq 1}$ be a stationary sequence of nonnegative random variables with finite mean and set $S_{n} = \sum_{i=1}^{n}T_{i}$. In \citet{Resnick}, the number of customers, $N_n$, in a system seen by the $n$th arriving customer is defined by $N_0 = 0$ and
\[
    N_{n} = ( N_{n-1}+1-\Gamma( S_{n-1},S_{n} ) )_{+}, \qquad n\geq 1,
\]
where $\Gamma$ is a homogeneous Poisson process with intensity $\mu$ independent of $\{T_{i}\}_{i\geq 1}$. The waiting time, $W_{n}$, of the $n$th arriving customer satisfies
\[
    W_{n} \overset{d}{=} \sum_{i=1}^{N_{n}+1}X_{i}
\]%
where $X_{1},X_{2},\ldots $ are iid exponentially distributed random variables with common mean $\mu$ independent of $N_{n}$. If $\{T_{i}\}_{i\geq 1}$ is a reversible, stationary, ergodic process and if $\E[T_{1}]<\mu $, then $N_{n}$ converges in distribution to a random variable $N$, the distribution of which may, under additional conditions, have a regularly varying tail. Whereas in \citet{Resnick} only a lower bound is given for the tail of
\[
    W \overset{d}{=} \sum_{i=1}^{N+1} X_{i},
\]%
Theorems~\ref{T:CV} and \ref{T:CV2} give sufficient conditions such that actually
\[
    \Pr [ W > x ] \sim \Pr [ N > \mu^{-1} x ], \qquad x \to \infty.
\]

\paragraph*{Teletraffic arrivals.}
A large number of teletraffic measurements shows that file sizes and transmission times exhibit heavy tails and long-range dependence. Standard models for explaining these empirically observed facts are the so-called ON/OFF model and the infinite source Poisson model. In \citet{FGMS}, a model is introduced that extends these standard models in a simple, but realistic
way. They assume that the first packet of a flow of data arrives at the point $\Gamma_{j}$ of a Poisson process with intensity $\lambda >0$. Flow $j$ then consists of $K_j$ packets, the $k$th of which arrives at time $Y_{jk} = \Gamma_{j} + S_{jk}$, where
\[
    S_{jk}=\sum_{i=1}^{k}X_{ji},\qquad 0\leq k\leq K_{j}.
\]
In \citet{FGMS}, it is assumed that the variables $X_{ji}$ form an
array of iid nonnegative random variables and the $K_{j}$ are iid
integer-valued random variables independent of the $X_{ji}$. Of
interest is the tail behavior of the total transmission time
$S_{jK_{j}}$ under the assumption that the tail of $X_{ji}$ or
$K_{j}$ is regularly varying. Proposition 4.3 of \citet{FGMS} gives
results similar as those in Theorem~\ref{T:CV2} above. In
\citet[Proposition~4.9]{FGMS}, the reverse problem is considered as
well: if the tail of $S_{jK_{j}}$ is regularly varying with index
$-\alpha$ for some $\alpha > 0$ and if the tail of $K_{j}$ is
heavier than that of $X_{ji}$, then what can be said about $K_{j}$?
For instance, if the tail of $S_{jK_{j}}$ is regularly varying with
index $-\alpha \neq -1$ and if $\Pr [ X_{1} > x] = o(\Pr[N>x])$ as
$x\to \infty$, then the tail of $K_{j}$ must be regularly varying
with index $-\alpha$ as well.


\end{document}